%% file: authorsample.tex
\documentclass[graybox]{svmult}

\usepackage{hyperref}
\usepackage{type1cm}        
\usepackage{makeidx}         
\usepackage{graphicx}        
\usepackage{multicol}        
\usepackage[bottom]{footmisc}\usepackage{caption}
\usepackage{algorithm,algpseudocode}
\usepackage{pgfplotstable}
\usepackage{pgfplots}
\usepackage{array,multirow,graphicx}
\usepackage{pgfplots}
    \usepgfplotslibrary{groupplots}
    \usetikzlibrary{matrix}
\DeclareMathAlphabet{\pazocal}{OMS}{zplm}{m}{n}    
\pgfplotsset{compat = 1.3}
\usepackage{wrapfig}

\usepackage{xcolor}
\definecolor{TR}{HTML}{6a4c93}
\definecolor{SN}{HTML}{c9182c}

\definecolor{TLTR0}{rgb}{0.63,0.12156,0.48236} 
\definecolor{TLTR1}{HTML}{ff924c}
\definecolor{TLTR2}{HTML}{ffd700}
\definecolor{TLTR3}{HTML}{8ac926}
\definecolor{TLTR4}{HTML}{1982c4}

\def \0{\mathbf{0}}			% 

\def \R{{\mathbb{R}}}			% real num
\def \N {\mathbb{N}}			% numbers 
\def \I{{\mathbf{I}}}			% subdivision matrix
\def \0{\vec{0}}			% 

\def \R{\mathbb{R}}			% real numbers
\def \N {\mathbb{N}}			% nat. numbers

\def \I{\mathbf{I}}			% subdivision matrix

\usepackage{newtxtext}       
\usepackage[varvw]{newtxmath}       

\algnewcommand\algorithmiconput{\textbf{Constants:}}
\algnewcommand\algorithmicinput{\textbf{Input:}}
\algnewcommand\algorithmicoutput{\textbf{Output:}}
\algnewcommand{\algorithmicgoto}{\textbf{go to}}%

\algnewcommand\Constants{\item[\algorithmiconput]}
\algnewcommand\Input{\item[\algorithmicinput]}%
\algnewcommand\Output{\item[\algorithmicoutput]}%
\algnewcommand{\Goto}[1]{\algorithmicgoto~\ref{#1}}%

\makeindex

\begin{document}

\title*{Two-level trust-region method with random subspaces}
\author{Andrea Angino\orcidID{0009-0000-8525-375X},\\
Alena Kopaničáková\orcidID{0000-0001-8388-5518}, and\\
Rolf Krause\orcidID{0000-0001-5408-5271}.}
\institute{Andrea Angino \at  Università della Svizzera italiana, 
UniDistance Suisse; \email{andrea.angino@unidistance.ch}
\and Alena Kopaničáková \at  Università della Svizzera italiana, Université de Toulouse, INP-ENSEEIHT, IRIT, ANITI;  \email{alena.kopanicakova@usi.ch} \\
\and Rolf Krause \at Università della Svizzera italiana,
UniDistance Suisse,
King Abdullah University of Science and Technology;
\email{rolf.krause@kaust.edu.sa}}

\maketitle

\abstract
{\input{chapters/abstract}}

%\abstract{TLTR
%\input{chapters/abstract}
%}

\input{chapters/intro}

%\newpage
\input{chapters/method}
\input{chapters/num_results}

\vspace{-0.2cm}

\begin{acknowledgement}
This work was partially supported by the Swiss Platform for Advanced Scientific Computing (PASC) project ExaTrain and by the Swiss National Science Foundation through the projects ``ML$^2$ - Multilevel and Domain Decomposition Methods for Machine Learning'' (197041). 
\end{acknowledgement}

 \end{document}

%% file: chapters/abstract.tex
%We introduce a novel variant of the Magical Trust Region (MTR) method designed for solving unconstrained non-linear optimization problems exploiting two distinct directions to enhance optimization efficiency. 
%The first direction is derived from a classical model-based approach, ensuring accurate and reliable step computations. 
%The second direction is obtained via a technique called sketching, which significantly reduces the dimensionality of the associated trust-region subproblem, thereby lowering computational costs. 
%Our method combines these two directions within the trust region framework to achieve robust and efficient convergence. 
%The integration of the sketched direction with the classical model-based direction enables a balanced approach that leverages the strengths of both techniques.
% A detailed discussion of the methodology, including the construction and integration of the sketching direction, is provided. 
% Additionally, we validate the effectiveness and robustness of the MTR method through numerical experiments on various network optimization problems.

We introduce a two-level trust-region method (TLTR) for solving unconstrained nonlinear optimization problems. 
Our method uses a composite iteration step, which is based on two distinct search directions. 
The first search direction is obtained through minimization in the full/high-resolution space, ensuring global convergence to a critical point. 
The second search direction is obtained through minimization in the randomly generated subspace, which, in turn, allows for convergence acceleration. 
The efficiency of the proposed TLTR method is demonstrated through numerical experiments in the field of machine learning.

%% file: chapters/intro.tex
\section{Introduction}
We  consider  minimization problems of the following type:
\begin{equation}\label{eq:Problem}    \min_{x\in\mathbb{R}^n}f(x),
\end{equation} 
where $f$ is the twice differentiable objective function (bounded from below) and $n$ possibly very large. 
In practice, optimization problems of this type are often solved using trust-
region (TR) methods~\cite{Conn TRmethods}, which provide
global convergence also in the presence of non-convexity. 
TR methods is generate search directions by 
building a local quadratic model of the objective function and minimizing this model while adhering to a step-size constraint at each iteration.
The acceleration of TR methods through  multilevel approaches (Recursive 
Multilevel Trust Region or RMTR) was condsidered in~\cite{RMTR, Gross2009, 
kopanivcakova2020recursive, kopanivcakova2019subdivision, kopanicakova2022globally, gratton2023multilevel}.
These methods aim to reduce the computational cost of TR methods while 
ensuring their algorithmic scalability.
This is achieved by exploring computationally cheaper, albeit less accurate, models of $f$ on "coarser" levels, in addition to the standard quadratic models used on the original "fine" level.

In this work, we introduce a 
new approach to construct the search directions in the TR algorithm, giving rise to a novel two-level TR (TLTR) method.
Methodologically, our TLTR method is built upon the "magic" TR 
framework~\cite[Section 10.4.1.]{Conn TRmethods}, which assumes the existence 
of an \emph{oracle/magic} that is used to improve the search directions 
obtained by minimizing the quadratic models.
This framework has been successfully explored in the literature in the context of constrained~\cite{Conn MTR} and non-smooth optimization problems~\cite{Fletcher1982, Yuan1985}.
Here, we extend it to unconstrained optimization problems by identifying the 
\emph{oracle/magic} steps with the \emph{coarse-level/subspace} steps, 
commonly utilized in multilevel/subspace optimization.
Thus, in each iteration, our TLTR method utilizes a composite search 
direction consisting of a step obtained by minimizing the quadratic model and a step obtained by minimizing a coarse-level/subspace model.
This approach conceptually differs from RMTR methods~\cite{RMTR}, which alternate between these two search directions instead of composing them.

From the application point of view, we aim to apply our TLTR method to the training of the machine learning models.
To this aim, we construct coarse-level/subspace models using randomly generated subspaces obtained through sketching~\cite{Lindenstrauss}.
It is worth noting that while sketching was initially proposed to reduce data dimensionality while controlling information loss, it has recently been also employed to reduce the computational cost of optimization methods. 
This has led to the development of several sketched first/second-order optimization methods, e.g., sketched gradient descent~\cite{sketchingGD}, Newton methods~\cite{newtonSketch, Cartis} or subspace methods~\cite{nick}.
To the best of our knowledge, a majority of the existing sketched optimization algorithms rely solely on sketched models throughout every iteration.
This imposes restrictions on the minimal size of the sketched subspace, as it directly affects the quality of the resulting search directions. 
In contrast, the TLTR method uses sketching only to improve upon full-space computations and therefore the size of sketched subspace does not impact the algorithm's convergence, only its speed.

This manuscript is organized as follows.  
In Section \ref{Sec:TLTR}, we introduce the proposed TLTR method. 
In Section \ref{Sec:NumericalEX}, we demonstrate the effectiveness of the TLTR method using numerical examples from the field of machine learning.

%% file: chapters/method.tex
\section{Two-level TR (TLTR) with random subspaces}
\label{Sec:TLTR}
In this section, we present the novel TLTR algorithm, which can be seen as the TR method that utilizes the composite search direction to find a solution of~\eqref{eq:Problem}. 
Motivated by the multigrid methods~\cite{trottenberg2000multigrid}, the composite direction is constructed by combining a step performed on the full space, called a smoothing step, and the step performed on the subspace called a coarse-level/subspace step.
%The proposed method can be see as a variant of magical TR algorithm presented in~\cite[Chapter 10]{Conn TRmethods}. 
%Our contribution is the identification of magical step with randomly generated coarse spaces and viewpoint as a subspace decomposition, in particular multilevel, method. 
%Moreover, to the best of our knowledge the use of sketching as an oracle and application to machine-learning applications have not been so far considered in the literature. 

\subsection{The TLTR algorithm}
Utilizing the magical TR framework~\cite{Conn TRmethods}, the TLTR method starts with an initial guess $x_0 \in \R^n$. 
At each $k^{\text{th}}$ iteration, the algorithm first performs a smoothing, achieved by approximately solving the following  quadratic minimization problem (QP): 
\begin{align}
& \min_{p_{k}^F \in \mathbb{R}^n} m_{k}^F \left( p_{k}^F\right) :=  f(x_k)+\left<\nabla f\left(x_k\right), p_{k}^F \right>+\frac{1}{2} \left< p_{k}^F, \nabla^2f(x_k) p_{k}^F \right>, \nonumber \\
 & \text{subject to} \quad \left\|p_{k}^F\right\| \leq \Delta_k,
 \label{eq:modelfull}
\end{align}
where the model~$m_{k}^F$ is constructed by utilizing a second order Taylor approximation of~$f$ around current iterate~$x_k$ and the symbol~$\Delta_k \in \mathbb{R}^+$ denotes a trust-region radius.
%. 
%Moreover, the symbol~$\Delta_k \in \mathbb{R}^+$ denotes a trust-region radius, which controls the size of~$p_{k}^F \in \R^n$.   

After the smoothing step is performed, the algorithm proceeds to obtain the subspace step.
To this aim, we consider an iteration-dependent restriction operator $S_k \in \mathbb{R}^{\ell \times n}$ and the prolongation operator $S_k^T  \in \mathbb{R}^{n \times \ell}$, where~$\ell \ll n$. 
These transfer operators are assembled using the sketching techniques~\cite{Lindenstrauss}.
% , and their role is to map quantities from the full space to a randomly sketched subspace of size~$\ell \ll n$ and vice-versa.
The algorithm then obtains the subspace step~$p^S_k \in \R^{\ell}$ by taking the advantage of a model~${m}_k^S$ of $f$, constructed around the iterate~${x_{k+1/2}:= x_{k} + p_{k}^F}$, i.e., the iterate obtained after the smoothing step.
Notably, the model ${m}^S_k$ can be constructed using various techniques, e.g., following $1^{\text{st}}/2^{\text{nd}}$-order consistency approaches outlined in~\cite{kopanicakova2022use}. 
For the purpose of this work, we utilize a simple restriction of the quadratic approximation of~$f$ to a subspace induced by $S_k$.
Thus, $p^S_k \in \R^{\ell}$ is obtained by solving the following sketched QP problem:
\begin{align}
& \min_{p^S_k \in \R^{\ell}} m_k^S\left( p^S_k\right) := f(x_{k+1/2}) + \left<S_k \nabla f(x_{k+1/2}), p^S_k \right>+ \frac{1}{2} \left<p^S_k, S_{k} \nabla^2 f(x_{k+1/2}) S_k^T p^S_k\right>, \nonumber \\
& \text{subject to}  \quad \left\|p^S_k \right\| \leq \Delta_k, \label{eq:modelcorse}
\end{align}
where the TR radius $\Delta_k$ again controls the size of~$p^S_k$.
This algorithmic design follows closely the RMTR algorithm~\cite{kopanicakova2022use}, where the size of the coarse-level correction is also controlled by means of $\Delta_k$. 

Subsequently, the subspace search direction~$p^S_k$ is transferred back to the full space. 
However, the prolongated search direction~$S_k^T p_k^S$ is not automatically considered by the TLTR algorithm, but only if it provides a decrease in the objective function, thus only if 
\begin{equation}
\label{eq:conditionstep}
f\left(x_k + p_k^F + \alpha_k S_k^T p_k^S\right) \leq    f\left(x_k + p_k^F\right). 
\end{equation}
Otherwise, the subspace search direction is discarded, i.e., we set $p^S_{k}$ to be the zero vector. 
Note, in order to improve quality of $S_k^T p_k^S$, one can utilize $\alpha_k \in \R^+$, obtained using a line-search strategy, such that $f(x_k + p_k^F + \alpha_k S_k^T p_k^S) < f(x_k + p_k^F + S_k^T p_k^S)$.

% Moreover, the quality of $S_k^T p_k^S$ can be further improved, such that $f(x_k + p_k^F + \alpha_k S_k^T p_k^S) < f(x_k + p_k^F + S_k^T p_k^S)$, where $\alpha_k \in \R^+$ is obtained using a line-search strategy.

Finally, the TLTR algorithm assesses the quality of the composite search direction $p_k := p_k^F+ \alpha_k S_k^T p^S_{k}$ by means of the TR ratio~$\varrho_k$, defined as
\begin{equation}
\label{eq:rho_tr}
\varrho_k := \frac{f\left(x_k\right)-f\left(x_k+p_k^F+ \alpha_k S_k^T p_{k}^S \right) }{m_k^F\left(x_k\right)-m_k^F\left(x_k+p_k^F\right)+f\left(x_k+p_k^F\right)-f\left(x_k+p_k^F+ \alpha_k S_k^Tp_{k}^S\right)}.
\end{equation}
Here, the nominator measures the objective function change before and after taking the composite step~$p_k$. 
The denominator measures the change in the model $m_k^F$ after the smoothing step~$p_k^F$ as well as the change in the objective function~$f$ after performing the subspace step.
If~$\varrho_k > \eta_1$, where $\eta_1 > 0$, then the composite step $p_k$ decreases the objective function sufficiently and it is safe to utilize it. 
Otherwise the composite search direction $p_k$ has to be discarded. 
In addition, the TR radius $\Delta_k$ has to be adjusted accordingly, i.e., enlarged, or shrunk depending on the value of $\varrho_k$.

We remark that if $p_{k}^S$ is a zero vector, the TLTR algorithm reduces to the standard, single-level, TR algorithm with 
$\varrho_k^{TR} := \frac{f\left(x_k\right)-f\left(x_k+p_k^F\right) }{m_k^F\left(x_k\right)-m_k^F\left(x_k+p_k^F\right)}$.
Moreover, if~$S_k p_{k}^S$ is a non-zero vector and the condition~\eqref{eq:conditionstep} holds, then $\varrho_k > \varrho^{TR}_k$, whenever $\varrho_k^{TR} < 1$. %~\cite{Conn TRmethods}
As a consequence, the composite step~$p_k$ might be accepted by the TLTR algorithm even if the standard TR algorithm would have rejected the full-space smoothing step $p_k^F$.
This algorithmic feature allows for the enhanced convergence of the TLTR algorithm compared to the standard TR method. 
Algorithm~\ref{ALG:TLTR} summarizes the described procedure.

\begin{algorithm}[t]
\footnotesize
\caption{\footnotesize Two-Level Trust Region (TLTR) Method}
\label{ALG:TLTR}
\begin{algorithmic}[1]
\Require {$ f: \R^n \rightarrow \R,  x_0 \in \R^n,  \Delta_{0} \in \R^+, \ell < n \in \N, k=0$}
\Constants { $0 < \eta_1 \leq \eta_2 < 1, \  0 < \gamma_1 \leq  \gamma_2 < 1$}
\While{not converged}
\State $p_k^F := \underset{\left\| p_k^F \right\| \leq \Delta_k}{\text{argmin}} m_k^F\left(p_k^F\right)$ \Comment{Obtain full-space search direction}
\State $x_{k+1/2} :=x_k + p_k^F$ 
\State
\State Construct $S_k$ via sketching 
\State $p_k^S :=  \underset{\| {p}_k^S \| \leq \Delta_k}{\text{argmin}} m_k^S\left({p}_k^S\right) $ \Comment{Obtain subspace search-direction }
\State 
\vspace{-0.4cm}
\Comment{Assess the quality of the subspace step}
\begin{align*}
& p_{k}^S = 
\begin{cases}
{p}_k^S, \quad &\text{if} \ f\left(x_k+p_k^F + \alpha_k S_k {p}_k^S\right) < f\left(x_k + p_k^F \right)    \\
0, \quad & \text{otherwise} \hspace{6.75cm} \textcolor{white}{.}
\end{cases}  &
\end{align*} %
% \Comment{Assess and improve the quality of the subspace step}
% \begin{align*}
% & p_{k}^S = 
% \begin{cases}
% \alpha_k{p}_k^S, \quad &\text{if} \ \left< S_kp_k^S,\nabla f(x_{k+1/2})\right> <0     \\
% 0, \quad & \text{otherwise} \hspace{6.75cm} \textcolor{white}{.}
% \end{cases}  &
% \end{align*}
\State Evaluate $\varrho_k$ as in~\eqref{eq:rho_tr} \Comment{Assess the quality of the composite trial step}
\vspace{-0.3cm}
\begin{align*}
 & x_{k+1} := 
\begin{cases}
x_{k} + p_k^F + \alpha_k S_k p_{k}^S, \quad &\text{if} \ \varrho_k > \eta_1   \\
x_{k}, \quad & \text{otherwise}
\end{cases} &
& \Delta_{k+1} := 
\begin{cases}
[\Delta_k, \infty), \quad &\text{if} \ \ \varrho_k \geq \eta_2   \\
[\gamma_2 \Delta_k, \Delta_k], \quad &\text{if} \ \ \varrho_k \in [\eta_1, \eta_2]   \\
[\gamma_1 \Delta_k, \gamma_2 \Delta_k], \quad &\text{if} \ \ \varrho_k < \eta_2  
\end{cases} 
\end{align*}
\State $k = k+1$
\EndWhile
\State \Return $x_{k+1}$
\end{algorithmic}
\end{algorithm}

\subsection{The computational cost of the TLTR method}
% \subsection{Implementation details}
\label{Sec:MagicSketch}
Compared to the standard TR algorithm, one iteration of the TLTR method is computationally more expensive. 
The increased computational cost is mostly associated with the random subspace minimization step. 
In particular, we sketch the subspaces by generating $S_k \in \R^{\ell \times n}$ using one of the following two sketching approaches:
\vspace{-0.1cm}
\begin{itemize}
\item\label{def:Gaussian} \textbf{Gaussian:} All entries of $S_k$ are independently distributed as $\mathcal{N}(0,\ell^{-1})$.

\item\label{def:shash} \textbf{S-hashing (sparse):} For each column $j \in \{1, \dots, n\}$, independently, we uniformly sample $s \ll \ell$ distinct row indices $i_1, \dots, i_s \in \{1, \dots, \ell\}$ without replacement. We then assign $(S_k)_{i_k, j} = \pm \frac{1}{\sqrt{s}}$ for $k = 1, \ldots, s$.
% Randomly assign $(S_k)_{i_k, j} = \pm \frac{1}{\sqrt{s}}$ for $k = 1, \ldots, s$.
\end{itemize}
\vspace{-0.1cm}
Gaussian approach generates a dense $S_k$, while s-hashing creates a sparse $S_k$ with $s$ non-zero entries per column. 
Moreover, s-hashing maintains the sparsity of the matrices it acts upon, allowing for efficient construction of~$m_k^S$ and the solution of~\eqref{eq:modelcorse} with a lower computational cost compared to the Gaussian approach.

Notably, obtaining~$p_k^S$ is significantly cheaper than $p_k^F$, as the sketched derivatives can be computed at reduced costs since only their subspace projections are required. 
Additionally, since~$\ell \ll n$, the cost of solving the subspace QP problem~\eqref{eq:modelcorse}  is significantly lower than~\eqref{eq:modelfull}.  
Here, we point out that, for global convergence of the TR/TLTR algorithm, it is sufficient to solve~\eqref{eq:modelfull} approximately, as long as resulting~$p_{k}^F$ fulfills the sufficient decrease condition, see~\cite{Conn TRmethods}.
 Consequently, in our numerical experiments, we solve~\eqref{eq:modelfull} using computationally inexpensive QP solvers.
% Specifically, we utilize the Cauchy Point (CP) method, which provides the minimizer of the model within the trust region along the projected steepest-descent direction. 
% Alternatively, we use a few iterations of the Steihaug-Toint conjugate gradient (ST-CG) method~\cite{Conn TRmethods}, designed to follow the piecewise linear path connecting the CG iterates, subject to the trust region constraint~\cite{TRsubproblems}.
% In contrast, there are no requirements on the necessary accuracy of solving~\eqref{eq:modelcorse}, as long as the resulting $S_k p_k^S$ satisfies~\eqref{eq:conditionstep}. 
Specifically, we utilize the Cauchy Point (CP) method, or a few iterations of the Steihaug-Toint conjugate gradient (ST-CG) method~\cite{Conn TRmethods}.
In contrast, there are no requirements on the necessary accuracy of solving~\eqref{eq:modelcorse}, as long as the resulting $S_k p_k^S$ satisfies~\eqref{eq:conditionstep}. 
However, due to the small dimensionality of the sketched problem, we opt to solve~\eqref{eq:modelcorse} as accurately as possible by the ST-CG method.
This algorithmic design closely mirrors the well-established multigrid approaches~\cite{trottenberg2000multigrid}, where smoothing is typically achieved using computationally inexpensive solution strategies that are particularly effective at reducing high-frequency components of the error. 
Moreover, it is standard practice to solve the subspace problems accurately due to their small size, which in turn allows the effective elimination of low-frequency components of the error.

% In conclusion, the overhead in the TLTR method, compared to the standard TR algorithm, primarily stems from the computational cost associated with the random subspace minimization step. Specifically, this cost includes \(\mathcal{O}(n^2 \ell)\) and \(\mathcal{O}(I \cdot \ell^2 + n \ell)\) , where \(n\) is the dimension of the original space, \(\ell\) is the reduced dimension, and \(I\) is the number of iterations required to solve the reduced problem.

% In summary, the TLTR method incurs additional overhead compared to the standard TR algorithm, mainly due to the random subspace minimization step. The associated costs are \(\mathcal{O}(n^2 \ell)\) and \(\mathcal{O}(I \cdot \ell^2 + n \ell)\), where \(n\) is the original dimension, \(\ell\) is the reduced dimension, and \(I\) is the number of iterations needed.

%% file: chapters/num_results.tex
\section{Numerical examples}
\label{Sec:NumericalEX}
We illustrate the numerical performance of the proposed TLTR method using classification problems.
 In particular, given a data set $D = {(z_i, y_i)}_{i=1}^N$, we consider
 \vspace{-0.2cm}
 \begin{itemize}
 \item \textbf{logistic loss}, defined as $f_{LL}(x) :=\sum_{i=1}^{N} \log \left( 1 + \mathrm{e}^{-y_i \left< x, z_i \right>} \right) + \frac{\lambda}{2} \|x\|^2$,  
 \item \textbf{least-square loss}, defined as
 $f_{LS}(x) := \frac{1}{N} \sum_{i=1}^N\left(y_i-\frac{\mathrm{e}^{\left< x, z_i \right>}}{1+ \mathrm{e}^{\left< x, z_i \right>}}\right)^2 + \frac{\lambda}{2} \|x\|^2$, 
 \end{itemize}
 where $\lambda=\frac{1}{N}$ represents the regularization term.
For both problems, we utilize Australian, Mushrooms, and Gisette datasets, from the LIBSVM database\footnote{https://www.csie.ntu.edu.tw/~cjlin/libsvmtools/datasets/binary.html}. Please, refer to Table~\ref{tab:dataset_characteristics} for the details regarding the number of samples (N), the parameter dimension (n), and the approximate condition number ($\kappa$) of the arising Hessians.

% \textcolor{red}{Add minimization problems type to description/captions}

    \begin{minipage}[c]{0.35\linewidth}
        \centering
        \footnotesize
        \captionof{table}{\footnotesize Properties of the datasets}
        \label{tab:dataset_characteristics}                
        \begin{tabular}{|c|c|c|c|c|}
            \hline
            \textbf{Dataset} & N & n & $\kappa_{LL}$& $\kappa_{LS}$ \\       
            \hline
            Australian & 621 & 14 & $10^6$& $10^4$  \\
            \hline
            Mushroom & 6,499 & 112 & $10^2$& $10^6$ \\
            \hline
            Gisette & 6,000 & 5,000 & $10^4$& $10^8$  \\
            \hline
        \end{tabular}
    \end{minipage}
\hfill
    \begin{minipage}[c]{0.57\linewidth}
  We first investigate the performance of the proposed TLTR method with respect to different sketching strategies.
For reference, the TLTR's performance is always compared to the standard TR method (without subspace step). 
For both methods, the full QP problems are solved using either  ST-CG or CP method. Moreover, we also 
%For the TLTR method, the sketched QP problems are always solved using the ST-CG method due to their small size. 
%Moreover, we also include a comparison with the sketched Newton (SN) method (Newton performed only on sketched subspace) with backtracking line-search, as proposed in~\cite[Algorithm 1]{Cartis}.
%In this case, the linear systems arising in each Newton iteration are solved using a direct solver, i.e., exactly. 
%For all numerical examples, we utilize random initial guesses and terminate the solution as soon as $\| \nabla f \| < 10^{-7}$.  
\end{minipage}        
\vspace{0.1cm}

\noindent %ST-CG or CP method.  
% Moreover, we also 
include a comparison with the sketched Newton (SN) method (Newton performed only on sketched subspace) with backtracking line-search, as proposed in~\cite[Algorithm 1]{Cartis}.
In this case, the linear systems arising in each Newton iteration are solved exactly.
For all experiments, we utilize random initial guesses and terminate the solution process as soon as $\| \nabla f \| < 10^{-7}$.

Figure~\ref{fig:T1bis} illustrates the behavior of the TLTR method for different subspace dimensions $\ell$ for the logistic loss $f_{LL}$ (first row) and the least square loss $f_{LS}$ (second row), with the Australian and Mushroom datasets. 
As expected, as the size of $\ell$ increases, the number of TLTR iterations decreases, but at the same time, the computational cost of solving the sketched QP subproblem increases. 
As a consequence, sketching subspace such that $\ell$ equals approximately to $20-30\%$ of $n$ strikes an optimal balance between the convergence speedup and the required computational cost. We also see that for least squares problem with the Mushroom dataset significantly small subspaces facilitate the fast convergence. 
Moreover, we can observe that by utilizing the CP method as the full-space QP solver, the performance of the TR solver rapidly deteriorates. 
In this case, the speedup obtained by the TLTR method is more prominent due to the use of the second-order information, albeit only within the sketched subspace.

We further investigate the impact of the parameter $s$, used within the s-hashing strategy, on the performance of the TLTR method. 
Figure~\ref{fig:T2} (left) illustrates the obtained results for $f_{LL}$ loss with the Mushroom dataset. 
As we can see, larger values of $s$, which give rise to a denser transfer operator, do not give rise to significantly faster convergence. 
Therefore, for the rest of the numerical experiments, the parameter $s$ is set to $\approx 10\%$ of the subspace parameter $\ell$.
% Figure~\ref{fig:T2} (left) illustrates the obtained results for a logistic regression example with the Mushroom dataset. 
% As we can see, larger values of $s$, which give rise to a denser transfer operator, do not give rise to significantly faster convergence. 
% Therefore, for the rest of the numerical experiments, the parameter $s$ is set to $\approx 10\%$ of the subspace parameter $\ell$.
Figure~\ref{fig:T2} (right) depicts the results for $f_{LS}$ loss for the TR and TLTR methods with the Gisette dataset over the fixed budget of $5,000$ iterations.
As we can observe, the TLTR method enables convergence to a more accurate solution of an order of magnitude.

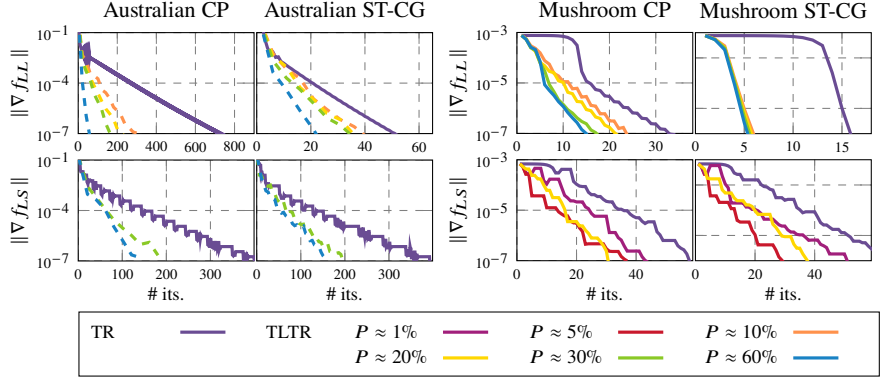
\begin{figure}
 \begin{tikzpicture}[]
    \begin{groupplot}[
        group style={
            % set how the plots should be organized
            group size = 5 by 2,
            % only show ticklabels and axis labels on the bottom
            % x descriptions at=edge bottom,
            % set the `vertical sep' to zero
            horizontal sep=1pt,
            vertical sep=0.35cm,
          },
	legend pos=north east,
        width=0.335\textwidth, 
        height=0.25\textwidth, 
	grid=major, % Display a grid
	grid style={dashed,gray}, % Set the style
	xmode=normal,
	ymode=log,
%	ylabel= $\pazocal{F}(\w) - \pazocal{F}(\w^*)$,
%	xmax= 35,
	ymin=1e-7, 
	ymax=1e-1,
	xmin=0,	
                ylabel shift = -0.1cm,
                xlabel shift = -0.1cm,
                title style={yshift = -0.12cm},         	
	tick label style={font=\tiny},
 	label style={font=\scriptsize},
	legend style={font=\footnotesize}        
      ]
      \nextgroupplot[align=left, 
      title={\footnotesize Australian CP},
      ylabel={$\| \nabla f_{LL}\|$}, 
      	xmax= 900]
%Iteration,GradNorm
       
\addplot[color = TR, very thick, solid] table[x=Iteration,y= GradNorm,col sep=comma] {results/Australian/Test4a/Line1.csv}; \label{pTR}  

\addplot[color = TLTR1, very thick, dashed, each nth point=20, filter discard warning=false, unbounded coords=discard ] table[x=Iteration,y= GradNorm,col sep=comma] {results/Australian/Test4a/Line2.csv}; 

\addplot[color = TLTR2, very thick, dashed, each nth point=20, filter discard warning=false, unbounded coords=discard] table[x=Iteration,y= GradNorm,col sep=comma] {results/Australian/Test4a/Line3.csv}; 
\addplot[color = TLTR3, very thick, dashed, each nth point=20, filter discard warning=false, unbounded coords=discard] table[x=Iteration,y= GradNorm,col sep=comma] {results/Australian/Test4a/Line4.csv}; 
\addplot[color = TLTR4, very thick, dashed, each nth point=20, filter discard warning=false, unbounded coords=discard] table[x=Iteration,y= GradNorm,col sep=comma] {results/Australian/Test4a/Line5.csv};

 \nextgroupplot[align=left, 
      title={\footnotesize Australian ST-CG},
      ylabel={}, 
        yticklabels={},      
      	xmax= 65]
\addplot[color = TR, very thick, solid] table[x=Iteration,y= GradNorm,col sep=comma] {results/Australian/Test4b/Line1.csv};

\addplot[color = TLTR1, very thick, dashed, each nth point=5, filter discard warning=false, unbounded coords=discard] table[x=Iteration,y= GradNorm,col sep=comma] {results/Australian/Test4b/Line2.csv}; 

\addplot[color = TLTR2, very thick, dashed, each nth point=5, filter discard warning=false, unbounded coords=discard] table[x=Iteration,y= GradNorm,col sep=comma] {results/Australian/Test4b/Line3.csv}; 
\addplot[color = TLTR3, very thick, dashed, each nth point=5, filter discard warning=false, unbounded coords=discard] table[x=Iteration,y= GradNorm,col sep=comma] {results/Australian/Test4b/Line4.csv}; 
\addplot[color = TLTR4, very thick, dashed, each nth point=5, filter discard warning=false, unbounded coords=discard] table[x=Iteration,y= GradNorm,col sep=comma] {results/Australian/Test4b/Line5.csv};

 \nextgroupplot[align=left, 
      title={\footnotesize },
        width=0.225\textwidth,
                axis lines=none, % Remove axis lines
                xtick=\empty, % Remove x-ticks
                ytick=\empty ]

      \nextgroupplot[align=left, 
      title={\footnotesize Mushroom CP},
	ymax=1e-3,      
      ylabel={$\| \nabla f_{LL}\|$}, 
        ytick={1e-3, 1e-5, 1e-7},      
      	xmax= 38]

\addplot[color = TR, very thick, solid] table[x=Iteration,y= GradNorm,col sep=comma] {results/mushroom/Test1/Line1.csv};

\addplot[color = TLTR1, very thick, solid] table[x=Iteration,y= GradNorm,col sep=comma] {results/mushroom/Test1/Line2.csv}; \label{p10TLTR}

\addplot[color = TLTR2, very thick, solid] table[x=Iteration,y= GradNorm,col sep=comma] {results/mushroom/Test1/Line3.csv}; \label{p20TLTR}
\addplot[color = TLTR3, very thick, solid] table[x=Iteration,y= GradNorm,col sep=comma] {results/mushroom/Test1/Line4.csv}; \label{p30TLTR}
\addplot[color = TLTR4, very thick, solid] table[x=Iteration,y= GradNorm,col sep=comma] {results/mushroom/Test1/Line5.csv};\label{p60TLTR}

      \nextgroupplot[align=left, 
      title={\footnotesize Mushroom ST-CG},
      xmax=18,	
	ymax=1e-3,            
      yticklabels={}]
      \addplot[color = TR, very thick, solid] table[x=Iteration,y= GradNorm,col sep=comma] {results/mushroom/Test1b/Line1.csv};

\addplot[color = TLTR1, very thick, solid] table[x=Iteration,y= GradNorm,col sep=comma] {results/mushroom/Test1b/Line2.csv}; 

\addplot[color = TLTR2, very thick, solid] table[x=Iteration,y= GradNorm,col sep=comma] {results/mushroom/Test1b/Line3.csv}; 
\addplot[color = TLTR3, very thick, solid] table[x=Iteration,y= GradNorm,col sep=comma] {results/mushroom/Test1b/Line4.csv}; 
\addplot[color = TLTR4, very thick, solid] table[x=Iteration,y= GradNorm,col sep=comma] {results/mushroom/Test1b/Line5.csv};

%Beginning of 2nd line!!!

      %1%

      \nextgroupplot[align=left, 
      title={},
      xlabel= $\#$ its., % Set the labels
      ylabel={$\| \nabla f_{LS}\|$}, 
      	xmax= 399]
%Iteration,GradNorm
       
\addplot[color = TR, very thick, solid] table[x=Iteration,y= GradNorm,col sep=comma] {results/LSreg/australianCP/Line1.csv};

\addplot[color = TLTR3, very thick, dashed, each nth point=20, filter discard warning=false, unbounded coords=discard] table[x=Iteration,y= GradNorm,col sep=comma] {results/LSreg/australianCP/Line2.csv}; 
\addplot[color = TLTR4, very thick, dashed, each nth point=20, filter discard warning=false, unbounded coords=discard] table[x=Iteration,y= GradNorm,col sep=comma] {results/LSreg/australianCP/Line3.csv};

%2%

 \nextgroupplot[align=left, 
      title={},
      ylabel={}, 
      xlabel= $\#$ its., % Set the labels
        yticklabels={},      
      	xmax= 399]
\addplot[color = TR, very thick, solid] table[x=Iteration,y= GradNorm,col sep=comma] {results/LSreg/australian/Line1.csv};

\addplot[color = TLTR3, very thick, dashed, each nth point=5, filter discard warning=false, unbounded coords=discard] table[x=Iteration,y= GradNorm,col sep=comma] {results/LSreg/australian/Line2.csv}; 
\addplot[color = TLTR4, very thick, dashed, each nth point=5, filter discard warning=false, unbounded coords=discard] table[x=Iteration,y= GradNorm,col sep=comma] {results/LSreg/australian/Line3.csv};

%3%
 \nextgroupplot[align=left, 
      title={\footnotesize },
        width=0.225\textwidth,
                axis lines=none, % Remove axis lines
                xtick=\empty, % Remove x-ticks
                ytick=\empty ]

      \nextgroupplot[align=left, 
      title={},
	ymax=1e-3,    
xlabel= $\#$ its., % Set the labels 
      ylabel={$\| \nabla f_{LS}\|$}, 
        ytick={1e-3, 1e-5, 1e-7},      
      	xmax= 59]

\addplot[color = TR, very thick, solid] table[x=Iteration,y= GradNorm,col sep=comma] {results/LSreg/mushroomCP/Line1.csv};
\addplot[color = TLTR0, very thick, solid] table[x=Iteration,y= GradNorm,col sep=comma] {results/LSreg/mushroomCP/Line1perc.csv};

\addplot[color = SN, very thick, solid] table[x=Iteration,y= GradNorm,col sep=comma] {results/LSreg/mushroomCP/Line2.csv}; \label{p5TLTR}

\addplot[color = TLTR2, very thick, solid] table[x=Iteration,y= GradNorm,col sep=comma] {results/LSreg/mushroomCP/Line3.csv};
%4%

      \nextgroupplot[align=left, 
      title={},
      xmax=59,	
	ymax=1e-3,            
 xlabel= $\#$ its., % Set the labels
      yticklabels={}]
      \addplot[color = TR, very thick, solid] table[x=Iteration,y= GradNorm,col sep=comma] {results/LSreg/mushroom/Line1.csv};
   \addplot[color = SN, very thick, solid] table[x=Iteration,y= GradNorm,col sep=comma] {results/LSreg/mushroom/Line5perc.csv};
\addplot[color = TLTR0, very thick, solid] table[x=Iteration,y= GradNorm,col sep=comma] {results/LSreg/mushroom/Line2.csv};\label{pgfplots:SN0}

\addplot[color = TLTR2, very thick, solid] table[x=Iteration,y= GradNorm,col sep=comma] {results/LSreg/mushroom/Line20.csv}; 

\end{groupplot}
    \matrix [ draw, matrix of nodes, anchor = west, node font=\scriptsize,
    column 1/.style={nodes={align=left,text width=1.0cm}},
    column 2/.style={nodes={align=left,text width=1.0cm}},    
    column 3/.style={nodes={align=left,text width=1.0cm}},
    column 4/.style={nodes={align=left,text width=1.0cm}},
    column 5/.style={nodes={align=left,text width=1.0cm}},
    column 6/.style={nodes={align=left,text width=1.0cm}},
    column 7/.style={nodes={align=left,text width=1.0cm}},
    column 8/.style={nodes={align=left,text width=1.0cm}},
    column 9/.style={nodes={align=left,text width=1.0cm}},
    column 10/.style={nodes={align=left,text width=1.0cm}},
    column 11/.style={nodes={align=left,text width=1.0cm}},  
    column 12/.style={nodes={align=left,text width=1.0cm}},
    column 13/.style={nodes={align=left,text width=1.0cm}},      
    ] at (0.0, -2.8)
    {
   TR &  \ref{pTR}  & & TLTR&    &  $P \approx 1\%$& \ref{pgfplots:SN0} &  $P \approx 5\%$& \ref{pgfplots:SN} & $P \approx 10\%$ &\ref{p10TLTR}   \\
   % TLTR    &  $P \approx 10\%$ &\ref{p10TLTR}  &
       &   &  & &  &  
   $P  \approx 20\%$ &\ref{p20TLTR}  &
    $P  \approx 30\%$ &\ref{p30TLTR}  &
   $P  \approx 60\%$ &\ref{p60TLTR}  \\  
    };

  \end{tikzpicture}
  \caption{ \footnotesize
  Convergence history of TR and TLTR for solving~\eqref{eq:Problem} with $f_{LL}$ (first row) and $f_{LS}$ (second row) with Gaussian (dashed lines, Australian dataset) and s-hashing (solid lines, Mushroom, $s = \left\lceil \ell/4\right\rceil$) for subspaces of varying size~$\ell = \lceil n P\rceil$, where $P$ denotes the portion of the full-space parameters.
To solve QP problems on full space, we use 2 iterations of ST-CG, or CP methods.}
  \label{fig:T1bis}
\end{figure}
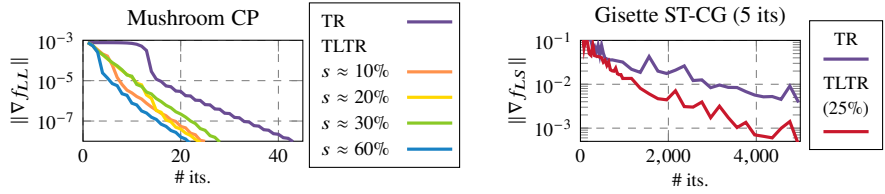
\begin{figure}
    \begin{tikzpicture}[]
            \matrix [ draw, matrix of nodes, anchor = west, node font=\scriptsize,
        column 1/.style={nodes={align=left,text width=0.95cm}}, % Reduced text width
        column 2/.style={nodes={align=center,text width=0.555cm}}, % Reduced text width
        ] at (3.0, 0.75) % Adjusted position for smaller legend
        {
        TR & \ref{pgfplot1:tr}  \\
        TLTR \\
        $s \approx 10\%$ &\ref{pgfplot1:TLTR1}  \\
        $s \approx 20\%$ &\ref{pgfplot1:TLTR2}  \\
        $s \approx 30\%$ &\ref{pgfplot1:TLTR3}  \\
        $s\approx 60\%$ &\ref{pgfplot1:TLTR4}  \\  
        };
    \matrix [ draw, matrix of nodes, anchor = west, node font=\scriptsize,
    column 0.8/.style={nodes={align=left,text width=1cm}}, % Adjusted text width for readability
    ] at (9.65, 0.75) % Adjusted position for the smaller, vertical legend
    {
    TR \\ \ref{pgfplot2:tr}  \\
    TLTR \\ (25\%)\\\ref{pgfplot2:TLTR1}  \\
    };        
        \begin{groupplot}[
            group style={
                group size=2 by 1,
                horizontal sep=105pt,
                vertical sep=0.5cm,
            },
            legend pos=south west, % Legend position
            width=0.385\textwidth, % Slightly reduced width
            height=0.25\textwidth, % Slightly reduced height
            grid=major,
            grid style={dashed,gray},
            xmode=normal,
            ymode=log,
            xlabel={$\#$ its.},
            ylabel={$\| \nabla f_{LL}\|$},
            ymin=1e-8, 
            ymax=10e-4,
            xmin=0,    
            ylabel shift = -0.1cm,
            xlabel shift = -0.1cm,
            title style={yshift = -0.12cm},                
            tick label style={font=\scriptsize},
            label style={font=\scriptsize},
            legend style={font=\scriptsize, inner sep=1pt, nodes={scale=0.8, transform shape}}, % Reduced legend size
        ]
            \nextgroupplot[
                align=left, 
                title={\footnotesize Mushroom CP}, % Title of the second plot
                xmax=45,
            ]
            \addplot[color=TR, very thick, solid] table[x=Iteration, y=GradNorm, col sep=comma] {results/mushroom/Test2/Line1.csv};
            \label{pgfplot1:tr}
            
            \addplot[color=TLTR1, very thick, solid] table[x=Iteration, y=GradNorm, col sep=comma] {results/mushroom/Test2/Line2.csv}; 
            \label{pgfplot1:TLTR1}
            
            \addplot[color=TLTR2, very thick, solid] table[x=Iteration, y=GradNorm, col sep=comma] {results/mushroom/Test2/Line3.csv}; 
            \label{pgfplot1:TLTR2}
            
            \addplot[color=TLTR3, very thick, solid] table[x=Iteration, y=GradNorm, col sep=comma] {results/mushroom/Test2/Line4.csv}; 
            \label{pgfplot1:TLTR3}
            
            \addplot[color=TLTR4, very thick, solid] table[x=Iteration, y=GradNorm, col sep=comma] {results/mushroom/Test2/Line5.csv};
            \label{pgfplot1:TLTR4}

        \nextgroupplot[
            align=left, 
            title={\footnotesize Gisette ST-CG (5 its)}, % Updated title to Gisette
            xmax=4999,
        xlabel={$\#$ its.},
        ylabel={$\| \nabla f_{LS}\|$},
        ymin=5e-4, 
        ymax=1e-1,
        xmin=0,    
        ylabel shift = -0.1cm,
        xlabel shift = -0.1cm,            
        ]
        \addplot[color=TR, very thick, solid] table[
            x expr={\thisrow{Iteration}},
            y expr={mod(\thisrow{Iteration},1) == 0 ? \thisrow{GradNorm} : NaN},
            col sep=comma
        ] {results/LSreg/gis/15/Line1.csv};
        \label{pgfplot2:tr}
        
        \addplot[color=SN, very thick, solid] table[
            x expr={\thisrow{Iteration}},
            y expr={mod(\thisrow{Iteration},1) == 0 ? \thisrow{GradNorm} : NaN},
            col sep=comma
        ] {results/LSreg/gis/15/Line2.csv}; 
        \label{pgfplot2:TLTR1}

    \end{groupplot}
\end{tikzpicture}
\caption{\footnotesize Left: Convergence of the TLTR method with $s$-hashing strategy ($\ell=\left\lceil n/5 \right\rceil$) for different values of sketching parameter $s$ for logistic loss with Mushroom dataset. 
Right: Convergence history of TR and TLTR for the least-square loss minimization problem with Gisette dataset.}
\label{fig:T2}
\end{figure}

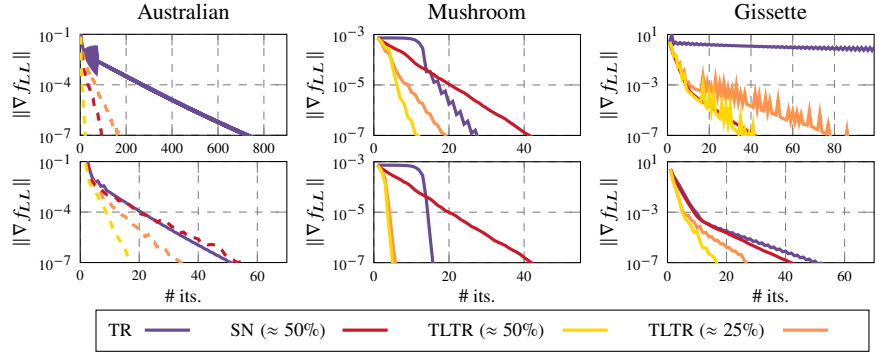
\begin{figure}
 \begin{tikzpicture}[]
    \begin{groupplot}[
        group style={
            % set how the plots should be organized
            group size = 5 by 2,
            % only show ticklabels and axis labels on the bottom
            % x descriptions at=edge bottom,
            % set the `vertical sep' to zero
            horizontal sep=1.5pt,
            vertical sep=0.35cm,
          },
	legend pos=north east,
        width=0.37\textwidth, 
        height=0.25\textwidth, 
	grid=major, % Display a grid
	grid style={dashed,gray}, % Set the style
	xmode=normal,
	ymode=log,
	% xlabel= $\#$ its., % Set the labels
%	ylabel= $\pazocal{F}(\w) - \pazocal{F}(\w^*)$,
%	xmax= 35,
	ymin=1e-10, 
	ymax=10,
	xmin=0,	
	tick label style={font=\tiny},
 	label style={font=\scriptsize},
	legend style={font=\footnotesize}, 
                ylabel shift = -0.1cm,
                xlabel shift = -0.1cm,
                title style={yshift = -0.12cm},         		
      ]
      \nextgroupplot[align=left, 
      title={\footnotesize Australian},
      ylabel={$\| \nabla f_{LL}\|$}, 
	ymin=1e-7, 
	ymax=1e-1,      
      	xmax= 900] %Iteration,GradNorm
\addplot[color = TR, very thick, solid] table[x=Iteration,y= GradNorm,col sep=comma] {results/Australian/Test6a/Line1.csv};\label{pgfplots:TR}

\addplot[color = TLTR1, very thick, dashed, each nth point=5, filter discard warning=false, unbounded coords=discard] table[x=Iteration,y= GradNorm,col sep=comma] {results/Australian/Test6a/Line2.csv}; 

\addplot[color = SN, very thick, dashed, each nth point=5, filter discard warning=false, unbounded coords=discard] table[x=Iteration,y= GradNorm,col sep=comma] {results/Australian/Test6a/Line3.csv}; 
\addplot[color = TLTR2, very thick, dashed, each nth point=5, filter discard warning=false, unbounded coords=discard] table[x=Iteration,y= GradNorm,col sep=comma] {results/Australian/Test6a/Line4.csv}; 
 \nextgroupplot[align=left, 
      title={\footnotesize },
        width=0.225\textwidth,
                axis lines=none, % Remove axis lines
                xtick=\empty, % Remove x-ticks
                ytick=\empty ]

      \nextgroupplot[align=left, 
      title={\footnotesize Mushroom},
	ylabel={$\| \nabla f_{LL}\|$},
        ytick={1e-3, 1e-5, 1e-7},     
	ymin=1e-7, 
	ymax=1e-3,                 
      	xmax= 55]
\addplot[color = TR, very thick, solid] table[x=Iteration,y= GradNorm,col sep=comma] {results/mushroom/Test3/Line1.csv};

\addplot[color = TLTR1, very thick, solid] table[x=Iteration,y= GradNorm,col sep=comma] {results/mushroom/Test3/Line2.csv}; 

\addplot[color = SN, very thick, solid] table[x=Iteration,y= GradNorm,col sep=comma] {results/mushroom/Test3/Line3.csv}; 
\addplot[color = TLTR2, very thick, solid] table[x=Iteration,y= GradNorm,col sep=comma] {results/mushroom/Test3/Line4.csv};
      
 \nextgroupplot[align=left, 
      title={\footnotesize },
        width=0.225\textwidth,
                axis lines=none, % Remove axis lines
                xtick=\empty, % Remove x-ticks
                ytick=\empty ]

 \nextgroupplot[align=left, 
      title={\footnotesize Gissette},
      xmax=99,	
	ymin=1e-7, 
	ylabel={$\| \nabla f_{LL}\|$},	
	ymax=10,            
      ]
     \addplot[color = TR, very thick, solid] table[x=Iteration,y= GradNorm,col sep=comma] {results/Gisette/Test9a/Line1.csv};

\addplot[color = TLTR1, very thick, solid] table[x=Iteration,y= GradNorm,col sep=comma] {results/Gisette/Test9a/Line2.csv}; 

\addplot[color = SN, very thick, solid] table[x=Iteration,y= GradNorm,col sep=comma] {results/Gisette/Test9a/Line3.csv}; %\label{pgfplots:SN}
\addplot[color = TLTR2, very thick, solid] table[x=Iteration,y= GradNorm,col sep=comma] {results/Gisette/Test9a/Line4.csv};

\nextgroupplot[align=left, 
      title={\footnotesize },
      xlabel= $\#$ its., % Set the labels
	ymin=1e-7, 
	ymax=1e-1,           
      ylabel={$\| \nabla f_{LL}\|$},  
      xmax= 70]
%Iteration,GradNorm
       
\addplot[color = TR, very thick, solid] table[x=Iteration,y= GradNorm,col sep=comma] {results/Australian/Test6b/Line1.csv};

\addplot[color = TLTR1, very thick, dashed] table[x=Iteration,y= GradNorm,col sep=comma] {results/Australian/Test6b/Line2.csv}; 

\addplot[color = SN, very thick, dashed] table[x=Iteration,y= GradNorm,col sep=comma] {results/Australian/Test6b/Line3.csv}; 
\addplot[color = TLTR2, very thick, dashed] table[x=Iteration,y= GradNorm,col sep=comma] {results/Australian/Test6b/Line4.csv};

 \nextgroupplot[align=left, 
 xlabel= $\#$ its., % Set the labels
      title={\footnotesize },
        width=0.225\textwidth,
                axis lines=none, % Remove axis lines
                xtick=\empty, % Remove x-ticks
                ytick=\empty ]

  \nextgroupplot[align=left, 
  xlabel= $\#$ its., % Set the labels
      title={\footnotesize },
	ylabel={$\| \nabla f_{LL}\|$},
        ytick={1e-3, 1e-5, 1e-7},     	
	ymin=1e-7, 
	ymax=1e-3,           
      	xmax= 55]
\addplot[color = TR, very thick, solid] table[x=Iteration,y= GradNorm,col sep=comma] {results/mushroom/Test3b/Line1.csv};

\addplot[color = TLTR1, very thick, solid] table[x=Iteration,y= GradNorm,col sep=comma] {results/mushroom/Test3b/Line2.csv}; 

\addplot[color = SN, very thick, solid] table[x=Iteration,y= GradNorm,col sep=comma] {results/mushroom/Test3b/Line3.csv}; 
\addplot[color = TLTR2, very thick, solid] table[x=Iteration,y= GradNorm,col sep=comma] {results/mushroom/Test3b/Line4.csv};

\nextgroupplot[align=left, 
      title={\footnotesize },
        width=0.225\textwidth,
                axis lines=none, % Remove axis lines
                xtick=\empty, % Remove x-ticks
                ytick=\empty ]

      \nextgroupplot[align=left, 
      title={\footnotesize },
      xlabel= $\#$ its., % Set the labels
      xmax=70,	
	ymin=1e-7, 
	ymax=10,   
	ylabel={$\| \nabla f_{LL}\|$}]
     \addplot[color = TR, very thick, solid] table[x=Iteration,y= GradNorm,col sep=comma] {results/Gisette/Test9/Line1.csv};

\addplot[color = TLTR1, very thick, solid] table[x=Iteration,y= GradNorm,col sep=comma] {results/Gisette/Test9/Line2.csv}; \label{pgfplots:TLTR}

\addplot[color = SN, very thick, solid] table[x=Iteration,y= GradNorm,col sep=comma] {results/Gisette/Test9/Line3.csv}; \label{pgfplots:SN}
\addplot[color = TLTR2, very thick, solid] table[x=Iteration,y= GradNorm,col sep=comma] {results/Gisette/Test9/Line4.csv};  \label{pgfplots:TLTRs}
\end{groupplot}

    \matrix [ draw, matrix of nodes, anchor = west, node font=\scriptsize,
    column 1/.style={nodes={align=left,text width=0.3cm}},
    column 2/.style={nodes={align=left,text width=1.0cm}},
    column 3/.style={nodes={align=left,text width=1.3cm}},
    column 4/.style={nodes={align=left,text width=1.0cm}},        
    column 5/.style={nodes={align=left,text width=1.6cm}},
    column 6/.style={nodes={align=left,text width=1.0cm}},        
    column 7/.style={nodes={align=left,text width=1.6cm}},
    column 8/.style={nodes={align=left,text width=1.0cm}},            
    ] at (0.2, -2.6)
    {
    TR & \ref{pgfplots:TR}  & 
    SN ($\approx 50\%$)& \ref{pgfplots:SN}  &     
    TLTR ($\approx 50\%$) & \ref{pgfplots:TLTRs} &
    TLTR ($\approx 25\%$) & \ref{pgfplots:TLTR}   \\
    };

  \end{tikzpicture}
  \caption{\footnotesize Convergence history of TR, SN and TLTR with Gaussian (dashed lines, Australian dataset) and s-hashing (solid lines, Gissette/Mushroom) subspaces.
The subspace sizes are chosen using a portion of full space $n$ stated in brackets. 
The QP problems on the full space are solved using CP/ST-CG (2/5 its) methods (Top/Bottom row).}
  \label{fig:T3s}
\end{figure}

%Figure \ref{fig:T3s} shows the results of Algorithm 1 in terms of the gradient norm reduction as a function of the number of iterations on the Australian, Mushroom, and Gisette datasets. The blue line represents the TR method, while the yellow one indicates the SN method. The remaining lines depict our TLTR strategy with different sketching configurations. Dashed lines correspond to Gaussian sketching (Definition 1), and solid lines use the $s$-hashing strategy (Definition 2). 
%All plot features different configurations, as detailed in Table~\ref{tab:line_color_datasets}

% AK:: read this again - does not read well
Figure~\ref{fig:T3s} illustrates the comparison of the TLTR method with respect to the TR and SN methods for logistic regression examples. 
Following the analysis presented in~\cite{Cartis}, the size of the sketched subspace for the SN method is set to be fairly large (50\% of $n$). For the TLTR method, we consider the same settings, as well as significantly smaller subspaces (25\% of $n$), guided by our previous experiments. As we can observe from the obtained results, the TLTR method outperforms the TR and SN methods. 
We attribute this behavior to taking advantage of both full-space and sketched search directions. 
In contrast, TR takes advantage of only full-space information, while SN utilizes only sketched quantities to obtain a search direction. 
Moreover, we also note that the observed speedups grow as the problems become larger and more ill-conditioned, i.e., as $n$ and $\kappa$ grow.
This highlights the potential of the proposed TLTR method for solving large-scale problems of practical importance, which we plan to explore in the future work.